\documentclass{amsart}
\usepackage{amssymb}
\usepackage{graphicx}
\usepackage{enumerate}
\newtheorem{theorem}{Theorem}[section]
\newtheorem{lemma}[theorem]{Lemma}
\newtheorem{prop}[theorem]{Proposition}

\theoremstyle{definition}

\theoremstyle{remark}
\newtheorem{remark}[theorem]{Remark}

\numberwithin{equation}{section}

\begin{document}

\title{$\alpha$-stable densities  are hyperbolically completely monotone for 
$\alpha\in ]0,1/4]\cup[1/3,1/2]$ }
\author{Sonia Fourati}
\address{INSA Rouen + LPMA\\
76131 St Etienne du Rouvray\\ FRANCE
}
\email{soniafourati@insa-rouen.fr}

%    General info
\date{04-09-2013}

\begin{abstract}
 We investigate the problem raised by L. Bondesson in \cite{bond},  about the hyperbolic complete monotonicity of $\alpha$-stable densities. We prove that densitites of subordinators of order $\alpha$ are HCM for $\alpha\in]0,1/4]\cup[1/3,1/2]$.
\end{abstract}
\maketitle
\section{Introduction} Hyperbolically completely monotone functions (HCM in short)
 were  introduced by L. Bondesson \cite{bond} in order to analyze infinitely divisible distributions, and particularly the so-called generalized gamma convolutions introduced by O. Thorin \cite{thor}. We recall their definition in section 1 below.

Bondesson showed that the densities of $\alpha$-stable positive random variables are HCM for $\alpha=n^{-1}$, for any integer $n\geq 2$. Furthermore, he conjectured that the HCM property actually holds for all $\alpha\in]0,1/2]$. Recently, Wissem Jedidi and
Thomas Simon \cite{WissSim}  investigated some aspects of the problem.   I thank them for  pointing out this question to me.

In this paper we prove this conjecture for values of $\alpha$ in $]0,1/4]\cup[1/3,1/2]$. 

For this we introduce the functions 
$$G_\alpha(x)=x^{-\frac{1}{\alpha}}g_\alpha(x^{-\frac{1-\alpha}{\alpha}})$$
where $g_\alpha$ is the density of the positive $\alpha$-stable distribution. We show that $G_\alpha$ extends to an analytic function on the slit plane
${\bf C}\setminus]-\infty,0]$. By analyzing its behaviour  at infinity and near the cut, we are able to prove that it has the following form
\begin{equation}\label{int-rep}
G_\alpha(z)=ce^{-\delta z}\exp\left(\int_0^{+\infty}\left[\frac{1}{z+t}-\frac{1}{1+t}\right]\theta(t)dt\right)
\end{equation}
where $c,\delta$ are  positive constants and $\theta$ takes values in $]0,1[$.

In order  that $G_\alpha$ be HCM it is then enough  that the function $\theta$ be increasing, which we prove for $\alpha\in[1/3,1/2]$. The HCM property for the remaining values of $\alpha$ is obtained by a multiplicative convolution argument.

This paper is organized as follows. In the section 2 we  recall some results of Zolotarev on densities of stable distributions. These are used in the next section to obtain the asymptotic behaviour of the function $G_\alpha$ in the complex plane. In section four we establish the integral representation (\ref{int-rep}). Finally, in section five, we prove that $\theta$ is increasing for $1/3<\alpha<1/2$, and we finish the proof of this part of the conjecture.

\section{Hyperbolically completely monotone functions}
We recall here the basic definition and properties of the class of  hyperbolically completely monotone functions, and refer to  \cite{bond} for more details.

  A real valued function $H$ defined on $]0,+\infty[$ is called hyperbolically completely monotone (HCM) if for every $u>0$ the function $H(uv)H(uv^{-1})$ is a completely monotone function of the variable $v+v^{-1}$. 
 Bondesson \cite{bond}  introduced this property in order to analyze infinitely divisible distributions, and particularly the so-called generalized gamma convolutions introduced by O. Thorin.

\begin{prop}\label{Proposition}\ \par
(i)  $H$ is HCM if and only if $H(x^{-1})$ is HCM.
\par
(ii) $H$ is HCM if and only if it admits the following representation
\begin{equation}\label{HCM}
H(x)=cx^{\beta-1}\exp\left(-a_1x-\int_1^\infty\log\frac{x+t}{1+t}\mu_1(dt)-a_2x^{-1}-
\int_1^\infty\log\frac{x^{-1}+t}{1+t}\mu_2(dt)\right)
\end{equation}
where $\beta$ is real, $a_1,a_2$ are positive constants and $\mu_1,\mu_2$ positive measures.\par
(iii) If $H$ is HCM then $H(x^{\beta})$ is HCM for all $\beta\leq 1$
\par
(iv)
$H$ is HCM if and only if the functions 
$x^{\gamma}H(x)$ are HCM for all values of $\gamma\in {\bf R}$.
\par
(v)
 If $X$ and $Y$ are independent positive random variables 
both with an HCM density then the random variable $XY$ also has an HCM density.
\end{prop}
In particular, from (iii) and (iv) we deduce that if $X$ is a positive random variable with HCM density, then $X^{\gamma}$ has HCM density for all $\gamma\geq 1$.
\section{Stable random variables}

Let $\alpha\in]0,1[$ and $\rho\in]0,1[$, we denote $g_{\alpha,\rho}$ the density of the strictly $\alpha$-stable distribution with asymmetry parameter $\rho$ (cf \cite{zol}). For $\rho=1$ (and only for this value) this distribution is supported on the half axis $]0,+\infty[$, and we 
simply put $g_\alpha=g_{\alpha,1}$.

The following result is an integral representation for the functions $g_{\alpha,\rho}$ on the positive axis, due to Zolotarev.

\begin{theorem}{(Zolotarev, \cite{zol}, Theorem 2.4.2)}

For all $x>0,\alpha,\rho\in]0,1[$
\begin{equation}\label{zolotarev}
g_{\alpha,\rho}(x)=(2i\pi)^{-1}\int_0^\infty
\frac{e^{-e^{-i\pi\rho\alpha}y^\alpha x^{-\alpha}}-e^{-e^{i\pi\rho\alpha}y^\alpha x^{-\alpha}}}{x}e^{-y}dy
\end{equation}
\end{theorem}

The following result which is easily obtained by a subordination argument,   plays an important role in the following.

\begin{lemma}
Let $X$ and  $Y$ be independent positive stable random variables, with respective parameters $(\alpha,\rho)$ and $(\beta,1)$, then
$XY^{1/\alpha}$ is a stable random variable with parameter $(\alpha\beta,\rho)$.
\end{lemma}

We deduce from the preceding lemma and Proposition \ref{Proposition} that 

\begin{prop}\label{Prop}
The set of $\alpha\in]0,1[$ such that $g_\alpha$ is HCM is a semigroup under multiplication.
\end{prop}
\section{The function $G_\alpha$}
 Denote $G_{\alpha}$ the function
\begin{equation}\label{Galpha}
G_{\alpha}(z)=(2i\pi)^{-1}\int_0^\infty \frac{e^{-e^{-i\pi\alpha}y^\alpha z^{1-\alpha}}-e^{-e^{i\pi\alpha}y^\alpha z^{1-\alpha}}}{z}e^{-y}dy\end{equation}
 where  we take  ( as  in the rest of the paper) for $z^h$,  the determination of the power function which is positive on $]0,+\infty[$ and analytic on ${\bf C}\setminus ]-\infty,0]$).

 This function $G_{\alpha}$
is analytic in ${\bf C}\setminus ]-\infty,0]$.
In fact
$z^\alpha G_\alpha(z)=F_\alpha(z^{1-\alpha})$ where $F_\alpha$ is an entire function. One has, for all $x>0$,
$$g_\alpha(x)=x^{-\frac{1}{1-\alpha}}G_\alpha(x^{-\frac{\alpha}{1-\alpha}})$$
and for all $z\in {\bf C}\setminus ]-\infty,0]$
$$G_\alpha(\bar z)=\overline {G_\alpha(z)}$$
For $r>0$ we denote
$$G_\alpha(-r^+)=\lim_{z\to-r,\Im(z)>0}G_\alpha(z)\qquad G_\alpha(-r^-)=\lim_{z\to-r,\Im(z)<0}G_\alpha(z)=\overline{G_\alpha(-r^+)}$$
the boundary values of $G_\alpha$.

\section{Behaviour near $0$.}
It follows from (\ref{Galpha}) that, as $z\to 0$,

\begin{equation}\label{G0}
G_\alpha(z)=\Gamma(\alpha+1)\frac{\sin(2\pi \alpha)}{\pi}z^{-\alpha}(1+O(|z|^{1-\alpha}))
\end{equation}

\section{Bounds at infinity}
\begin{theorem}\label{asymptotic}

Let $\theta\in ]-1,1[$ be fixed, and
$$\delta=(1-\alpha)\alpha^{\frac{\alpha}{1-\alpha}}\qquad c=(1-\alpha)^{-\frac{1}{2}}\alpha^{\frac{1}{2(1-\alpha)}}$$
then, as $r\to+\infty$, for $z=re^{i\pi\theta}$, one has
\begin{equation}\label{Galphaz}G_\alpha(z)\sim cz^{-\frac{1}{2}}e^{-\delta z}\end{equation}
As $r\to+\infty$
\begin{equation}\label{Galphar}G_\alpha(-r^+)\sim -icr^{-\frac{1}{2}}e^{\delta r}
\qquad
G_\alpha(-r^-)\sim icr^{-\frac{1}{2}}e^{\delta r}
\end{equation}

Furthermore, for some $R>0$, the function
$G_\alpha(z)z^{1/2}e^{\delta z}$
is uniformly bounded on ${\bf C}\setminus ]-\infty,0]\cap \{|z|>R\}$.
\end{theorem}

In order to obtain this asymptotic result, observe that one can rewrite the integral defining $G_\alpha$ as a contour integral:

$$G_\alpha(z)=(2i\pi)^{-1}\int_\Gamma \frac{e^{y-y^\alpha z^{1-\alpha}}}{z}dy$$
where $\Gamma$ is a contour which starts from $-\infty$, following the negative axis, taking the lower branch of $y^\alpha$, encircles $0$ then goes back to $-\infty$ along the negative axis, this time picking up the upper branch of $y^\alpha$.

In order to obtain the asymptotics we take $z=re^{i\theta}$ and rewrite the integral as
$$G_\alpha(z)=(2i\pi)^{-1}\int_\Gamma e^{r(y-y^\alpha e^{i\pi(1-\alpha)\theta})}e^{-i\pi\theta}dy$$
This integral is subject to the steepest descent method  (see \cite{Miller} for example) using the  unique saddle point at
$y=\alpha^{\frac{1}{1-\alpha}}e^{i\pi\theta}$ of the function $y-y^\alpha e^{i\pi(1-\alpha)\theta}$. This gives the point wise convergence  for a fixed $\theta$. In order to obtain the uniform convergence, first notice that uniform property is clear    for  $\theta$ in any compact subset of $]-1,+1[$, say for $\theta\in [-7/8, 7/8]$. then,  for $\theta \in ]7/8,1[$, the saddle point is over  the half line $]-\infty,0[$ and close to it, then one can  use another determination of $y^\alpha$ with a cut say on the half-line $\arg(y)=-3\pi/4$, and a contour encircling the cut and going back to a neighborhood of $-\infty$  by an arc with a ray going to infinity. Then  again one can  deform this contour to go through the saddle point and then conclude of the uniform convergence for 
$\theta \in ]7/8,1[$. A symmetrical argument gives the uniformity for $\theta \in ]-1, -{7\over 8}[$.

\section{Behaviour of $G_\alpha$ on the cut}

\begin{lemma}
For any $r>0$
\begin{equation}\label{r+}
G_\alpha(-r^+)=(2i\pi)^{-1}\int_0^\infty
\frac{e^{r^{1-\alpha}y^\alpha}-e^{e^{-2i\pi\alpha}r^{1-\alpha}y^\alpha}}{r}e^{-y}dy\end{equation}

\begin{equation}\label{rexp}
G_\alpha(-r^+)=(2i\pi)^{-1}\sum_1^\infty \frac{\Gamma(n\alpha +1)}{\Gamma(n+1)}(1-e^{-2i\pi n\alpha})r^{n(1-\alpha)-1}\end{equation}
\end{lemma}
{\sl Proof.} The first   formula follows at once  from (\ref{Galpha}) by letting $z\to -r$, the second one comes from expanding the exponentials in the numerator of (\ref{Galpha}) and integrating term by term. 
\qed

\begin{lemma}
For any $r>0$
one has $\Im(G_\alpha(-r^+))<0$. Furthermore, $-r^{\alpha}\Im(G_\alpha(-r^+))$ is an increasing function of $r$.
\end{lemma}

{\sl Proof.} By (\ref {rexp}) we get
$$-\Im(G_\alpha(-r^+))=(2\pi)^{-1}\sum_1^\infty \frac{\Gamma(n\alpha +1)}{\Gamma(n+1)}(1-\cos(2\pi n\alpha))r^{n(1-\alpha)-1}$$
in which all terms in the sum are positive;  the two claims are clear. 
\qed

In the sequel, denote 
 $$G_\alpha(-r^+)=: R(r)e^{-i\pi\theta(r)}$$ the polar decomposition of $G_\alpha(-r^+)$. Since $\Im(G_\alpha(-r^+))<0$, 
 one can chose $\theta(r)$ in $]0,1[$.

Observe also that $\theta(r)\to 1/2$ as $r\to+\infty$ (by (\ref{Galphar})).
 
\begin{remark} in fact one could also obtain from the integral representation that $\theta(t)-1/2=o(e^{-\epsilon r})$ as $r\to+\infty$ for some $\epsilon >0$, but we will not use this).
\end{remark}

\section{Integral representation}
\begin{prop}
For all $z\in{\bf C}\setminus ]-\infty,0]$
\begin{equation}
G_\alpha(z)=ae^{-\delta z}\exp\int_0^\infty \left[\frac{1}{z+t}-\frac{1}{1+t}\right]\theta(t)dt
\end{equation}
for some $a>0$.
\end{prop}

{\sl Proof.}  Let
$$L_\alpha(z)=\exp\int_0^\infty \left[\frac{1}{z+t}-\frac{1}{1+t}\right]\theta(t)dt$$

This is an analytic function on ${\bf C}\setminus ]-\infty,0]$, and it satifies, by well known properties of Stieltjes transforms,
$$\frac{L_\alpha(-r^+)}{L_\alpha(-r^-)}=e^{-2i\pi\theta(r)}$$
Furthermore, as $z\to\infty$, since $\theta(t)\to_{t\to+\infty} 1/2$, one has 
$$L_\alpha(z)=z^{-1/2}\exp(o(\log(|z|))$$
Near zero, one has 
$\theta(t)=\alpha+O(t^{1-\alpha})$ by (\ref{G0}), which implies
$$L_\alpha(z)\sim z^\alpha\qquad z\to 0$$
On the other hand, for $r>0$,
$$\frac{G_\alpha(-r^+)}{G_\alpha(-r^-)}=e^{-2i\pi\theta(r)}$$
therefore the function
$$E_\alpha(z)=e^{\delta z}G_\alpha(z)/L_\alpha(z)$$
is analytic on ${\bf C}\setminus ]-\infty,0]$, and can be extended continuouly to $\bf C\setminus\{0\}$. Since it is bounded near 0 it can be extended to  
 an entire function and it satifies 
$$E_\alpha(z)=\exp(o(\log(|z|))$$ at infinity thus it is  constant. Since both functions $G_\alpha, L_\alpha$ take positive values on $]0,+\infty[$, this constant is positive.\qed

\section{The function $\theta$ is monotone for $\alpha\in[1/3,1/2]$}
\begin{lemma} For $0\leq \rho\leq\inf(1,\frac{1}{2\alpha})$ the function
$\tilde g_{\alpha,\rho}(x)=x^{-1-\alpha}g_{\alpha,\rho}(x^{-1})$ is decreasing on $]0,+\infty[$. 
\end{lemma}
{\sl Proof.} Recall that if $X$ is a   stable variable  with parameters $(2\alpha,\rho)$, and $Y$    an independent stable variable with parameters $(1/2,1)$, then $Z=XY^{\frac{1}{2\alpha}}$ is a stable variable with parameters $(\alpha,\rho)$ . Since the density of $Y$ is 
$\frac{e^{-\frac{1}{2t}}}{\sqrt{2\pi t^3}}$ one has
$$g_{\alpha,\rho}(x)=2\alpha\int_0^\infty 
g_{2\alpha,\rho}(y)
\frac
{e^{-\frac{1}{2}(y/x)^{2\alpha}}y^{\alpha}}
{\sqrt{2\pi}x^{\alpha+1}}
dy$$ Therefore
$$x^{-1-\alpha}g_{\alpha,\rho}(x^{-1})=2\alpha\int_0^\infty g_{2\alpha,\rho}(y)
\frac{e^{-\frac{1}{2}(yx)^{2\alpha}}y^{\alpha}}{\sqrt{2\pi}}dy$$
which is clearly decreasing in $x$.

 \qed

\begin{lemma} For $\alpha\in[1/3,1/2]$  the function 
$r^\alpha\Re G_\alpha(-r^+)$ is decreasing.
\end{lemma}
{\sl Proof.} Note that, by formulas (\ref{zolotarev}) and (\ref{r+}) one has
$$\Re G_\alpha(-r^+)=r^{-1/\alpha}g_{\alpha,{1\over \alpha}-2}(r^{-\frac{1-\alpha}{\alpha}})$$
for $\alpha\in[1/3,1/2]$.
it follows that 
$$r^{\alpha}\Re G_\alpha(-r^+)=r^{\alpha-1/\alpha}g_{\alpha,{1\over \alpha} -2}(r^{\frac{1-\alpha}{\alpha}})=
x^{-1-\alpha}g_{\alpha,{1\over \alpha}-2}(x^{-1})$$ with $x=r^{1-\frac{1}{\alpha}}$.
The result follows from the preceding lemma.

\qed

\begin{theorem} 
For $\alpha\in]1/3,1/2|$, the function $\theta$ increases from the value $\alpha$ to the value $1/2$, and
$$G_\alpha(z)=\Gamma(\alpha+1)e^{-\delta z}z^{-1/2}\exp-\int_0^\infty \log(1+t/z)\theta'(t)dt$$

\end{theorem}

{\sl Proof.} One has
$$\tan(\pi\theta(r))=\frac{-r^\alpha\Im(G_\alpha(-r^+))}{r^\alpha\Re(G_\alpha(-r^+))}$$
and the numerator and denominator of this formula are positive and respectively increasing and decreasing. This implies that $\theta$ is increasing.
The other claim follows by integrating by parts. \qed

\section{The HCM property of stable distribution}

For $\alpha\in [\frac{1}{3},\frac{1}{2}]$ one has $\frac{1-\alpha}{\alpha}\geq 1$ and $G_\alpha$ is HCM. This implies that $g_\alpha$ is HCM. By Proposition (\ref{Prop}) the set of $\alpha$ such that $g_\alpha$ is HCM thus contains the multiplicative  semigroup generated by $[1/3,1/2]$, which is 
$]0,1/4]\cup [1/3,1/2]$.
\qed

\medskip

\begin{remark} Following the same arguments than above, one can prove  that for $\alpha\geq 1/2$, the function $\theta$ decreases from $\alpha$ to $1/2$ and consequently  $G_{\alpha}$ enjoys the next decomposition :

$$G_\alpha(z)=\Gamma(\alpha+1)e^{-\delta z}z^{-1/2}\exp \int_0^\infty \log(1+t/z)|\theta'(t)|dt$$

In other words, $e^{-\delta z} {1\over G_\alpha(z)}$ is an  HCM function.
\end{remark} 

%\tableofcontents


\begin{thebibliography}{10}

\bibitem{bond}
Bondesson, Lennart {\sl 
Generalized gamma convolutions and related classes of distributions and densities.} 
Lecture Notes in Statistics, 76. Springer-Verlag, New York, 1992. 
\bibitem{WissSim}
Jedidi, Wissem and  Simon Thomas {\sl Further examples of GGC and HCM functions}. To appear in Bernoulli Journal.
\bibitem {Miller} Miller Peter David {\sl Applied Asymptotic  Analysis}. AMS, 2006, Vol.77
 \bibitem{thor}
Thorin, Olof {\sl
On the infinite divisibility of the Pareto distribution.}
Scand. Actuar. J. 1977, no. 1, 31–40. 
\bibitem{zol}
Zolotarev, V. M. {\sl
One-dimensional stable distributions.}
Translated from the Russian by H. H. McFaden. Translation edited by Ben Silver. Translations of Mathematical Monographs, 65. American Mathematical Society, Providence, RI, 1986. 
   
\end{thebibliography}
\end{document}